# A note on the automorphism group of the root lattice of the U-dual modular group

Marcel Steiner[1]

*Abstract.* We study the inclusion system of the quantum deformed 2 dimensional Yang-Mills root module to the graded root module of the U-dual modular group. The irreducible representation of the U-dual modular group is the quantum deformed black brane throat. We expose the isomorphism between coherent superposition of the topological amplitude and the space of the automorphic forms of the U-dual modular group.

# 1. Introduction: From the index of Dirac operator to the union of the nilpotent orbits

We follow Atiah assumption based on the bijective homomorphism between a group variety and a vacuum variety. The congruence between number theory and homological features of the vacuum condensate is supporting for this reasoning.

There is a number atoms relation with a vacuum atoms such as fundamental correspondence between Riemann zeta function and the generating function of the black brane formulated by matrix theory where the eigenvalues $\lambda$ of the N x N matrix $\phi$ are just the roots of the zeta function. It is well known comprehend that Atiah-Bott-Shapiro isomorphism relates Clifford algebra to K-theory and that Bott periodicity theorem, which determines fundamental structure for K-theoretic group, is substantial for Atiah-Singer theorem and Thom isomorphism.

Naturally there is implication for inevitabillity of Langlands correspondence [GW, KW] and for the isomorphism between the cascading conifold throat and affine Weyl reflections of the ADE root system [CFKIV, FHHV, H]. We define a potential V for the ADE quiver variety by k-matrix models accordance with [DV, HSV]

$$Z_{k-cut} = \int \prod_k d\phi_k \ e^{-\Sigma_k \operatorname{Tr} V(\phi_k) \mathbf{g}^{-1}} (\prod_k \operatorname{Vol} \ (U(N)))^{-1}$$

or accordance with [DV2]. If Lie algebra has rank k, we have k matrices $\phi_a$ of rank $N_a$, one for each node of Dynkin diagram. Each matrix $\phi_a$ has an individual potential $V_a$. We define the generating function in terms of the eigenvalues $\lambda_{a,i}$ where $a = 1,...,k$ and $i = 1,...,N_a$

$$Z_{k-cut} = \int \prod_{a,i} d\lambda_{a,i} \ e^{\Sigma_{a,i} V_a(\lambda_{a,i}) \mathbf{g}^{-1}} \prod_{(a,i) \neq (b,j)} (\lambda_{a,i} - \lambda_{b,j})^{\alpha_a \times \alpha_b}$$

where $\alpha_a$ are the simple roots and $\alpha_a \times \alpha_b$ is Cartan matrix.

Let $S_{bh}$ is the black hole entropy and $C_k$ is $k^{th}$ Catalan number which counts the number of distinct ways the black hole centers can be generated

$$Z_{k-cut} = \sum_{k-1}^{\infty} C_{k-1} \sum_{p_k, q_k} e^{S_{bh}(p_1, q_1) + \cdots + S_{bh}(p_k, q_k)}$$

Our reasoning on the duality between a group variety and a vacuum variety is supported by the correspondence between the generalized Catalan numbers and Weyl group generators. We survey stringy phase transitions in the moduli space by the mirror to the topological vertex accordance with [BKMP] and we find exact agreement with the space of the automorphic forms determining the corresponding elliptic genus. The scheme of our observation is following. We deform the conifold by the k-flux-warped cascading throat and we study these transitions by Weyl group of the corresponding root system. We generalize conifold transitions by the condensation of Hagedorn spectrum of the dilaton condensate and we concentrate on the automorphism group of the appropriate root lattice. We generalize holographic description of k-center black hole throat [DGOV] and the quantum entanglement of the entropy of the k-center black hole [AOO].

In the line of our reasoning, the vacuum irreducible representation is just the quantum entanglement of the K-theoretic condensate determined by Bott periodicity theorem

$$\pi_n(O(\infty)) = \pi_{n+8}(O(\infty))$$

The base space of super-Yang-Mills/M-theory (and the appropriate vertex algebra) substantiates that matrix theory of universal SYM/M-theoretic backgrounds inherently implies nonassociative deformation of the brane worlvolume algebra: octonionic M2-, M5-instanton sectors are identical on nonassociative deformed algebraic variety.

Our research is focused on the dilaton groundstate orbits under the action of the full automorphism group of the root lattice of *the U-dual modular group* **M**. We postulate **M** as the isomorphism group generating total module of the vacuum torsion.

*The homomorphism between the root module of **M** and the superposition of the topological amplitude*

*and*

*the isomorphism between Weyl reflections of the root space of **M** and orbifold group variety of the generalized conifold*

is important proof for our postulate.

We construct the moduli space of the U-dual vacuum torsion as a one U-fold which is the union of the nilpotent orbits of **M**. The symmetry of the quantum entanglement of the dilaton condensate (which is the holographic dual of the only one unstable U-dual black brane throat) imposes background independent formulation of topological M-theory.

## 2. U-dual root module

Elementary observation is that the quantum deformed root system of a algebraic variety/a vacuum variety[(2)] is the function of the coupling costant **g**. We base on the algebraic variety of the quantum deformed 2 dimensional Yang-Mills/Chern-Simons theory ($qYM_2/qCS_2$) determining extremal black hole microstates [AOSV, BP].

**Definition 2.1** *Let $Z_{top}(U_1,...,U_{|2g-2|})$ is the topological string amplitude with $|2g-2|$ stacks of D-branes and $U_i$ is the holonomy group of the gauge field of $i^{th}$ stacks of D-branes*

$$Z_{bh} = \int \prod_{i=1}^{|2g-2|} dU_i \; |Z_{top}(U_1, \ldots, U_{|2g-2|})|^2$$

*Let $R_i$ are $A_\infty$ group representations and $Tr_{R_i} U_i$ (i=1,...,|2g-g|) are the characters of $R_i$*

$$Z_{top}(U_i) = \sum_{R_i} Z_{top\, R_i} \; Tr_{R_i} U_i$$

We explore that the space of the automorphic forms generating the partition function of qYM/qCS on Riemann surface in the holonomy basis

$$Z_{qYM/qCS} = \sum_R \dim_q(R) \; Tr_R U = \sum_R \chi_R(\mathbf{g}\rho) \; Tr_R U$$

(where $\dim_q(R)$ is the quantum deformed dimension of the group representation R, $\chi_R$ is Euler character of R and $\rho$ is Weyl vector) has the structure of a blossom flower in the moduli space of the vacuum dualities. Initial "bud-like" superset in the moduli space of the vacuum dualities is the space of the automorphic forms which generates $qYM_3/qCS_3$ modular form determining the generating function of topological M-theory. We base on the corresponding root lattice-valued basis vectors analysis and we discover the isomorphism between $qYM_3/qCS_3$ root module and its $qYM_2/qCS_2$ root submodule and

$$Z_{top} = Z_{top\, membrane} = Z_{top\, 3-brane} = \cdots$$

The inclusion system of the quantum deformed root submodules and generalized Freudenthal-Tits square modulo Bott periodicity theorem of the quantum deformed homotopy groups exposes the vacuum torsion module

|  | qℝ-valued base | qℂ -valued base | qℍ-valued base | q𝕆-value base |
|---|---|---|---|---|
| qℝ-valued base | q$\mathbf{A_1}$ root module | q$\mathbf{A_2}$ root module | q$\mathbf{C_3}$ root module | q$\mathbf{F_4}$ root module |
| qℂ -valued base | q$\mathbf{A_2}$ root module | q$\mathbf{A_2}$ ⊕ q$\mathbf{A_2}$ root module | q$\mathbf{A_5}$ root module | q$\mathbf{E_6}$ root module |
| qℍ-valued base | q$\mathbf{C_3}$ root module | q$\mathbf{A_5}$ root module | q$\mathbf{B_6}$ root module | q$\mathbf{E_7}$ root module |
| q𝕆-valued base | q$\mathbf{F_4}$ root module | q$\mathbf{E_6}$ root module | q$\mathbf{E_7}$ root module | q$\mathbf{E_8}$ root module |

$$\text{modulo q}\mathbb{O}: \pi_n(q\mathbf{D}_\infty) = \pi_{n+8}(q\mathbf{D}_\infty)$$

For importance of Bott periodicity theorem of the quantum deformed homotopy groups we analyze q$\mathbf{D_4}$ module homomorphism which is the product in the quantum deformed octonionic space q𝕆 [BI, HHM]. The automorphism group of q𝕆 is the quantum deformed group q$\mathbf{G_2}$

$$\text{Aut}(q\mathbb{O}) = q\mathbf{G}_2$$

which is the target space of the quantum deformed background of M-theory.

As a verification of a justifiability of our thoughts on the pure Platonism we build up the representation of the quantum deformed Clifford algebra qCliff$_8$, which classifies the quantum deformed charges of the vacuum fluxes, on the direct sum of 2 q$\mathbb{O}^{(3)}$, we find out the automorphism group of the quantum deformed exceptional Jordan algebra $\mathfrak{h}_3(q\mathbb{O})$, which counts particular class of the black brane microstates, as the quantum deformed group q$\mathbf{F_4}$

$$\text{Aut}(\mathfrak{h}_3(q\mathbb{O})) = q\mathbf{F}_4$$

and we construct the quantum deformed group q$\mathbf{E_8}$, on which self-dual root lattice live the basis vectors, by q𝕆-valued root module

$$\text{Tri}(q\mathbb{O}) \oplus \text{Tri}(q\mathbb{O}) \oplus (q\mathbb{O} \otimes q\mathbb{O})^3 = q\mathbf{E}_8$$

(where Tri(q𝕆) is the triality group of q𝕆).

**Definition 2.2** *There is the bijective homomorphism between the inclusion system of the submodules of the graded module of the vacuum torsion*

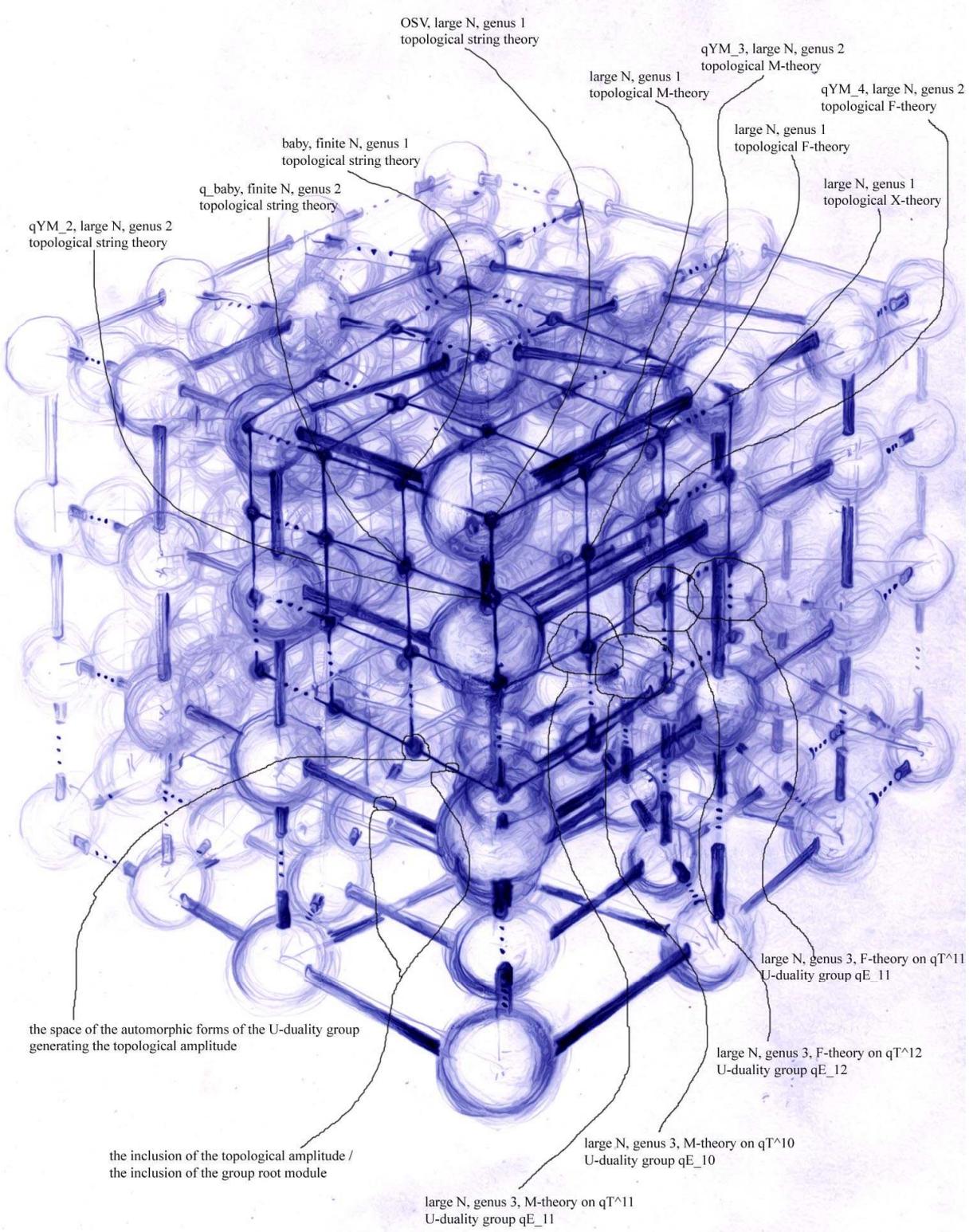

(where OSV if for the OSV conjecture, baby for the baby universes, $_q$baby for the quantum entanglement of the baby universes and topological X-theory for the supercritical topological M-theory)

*and the isomorphism web of the root system of the T-dual modular group*

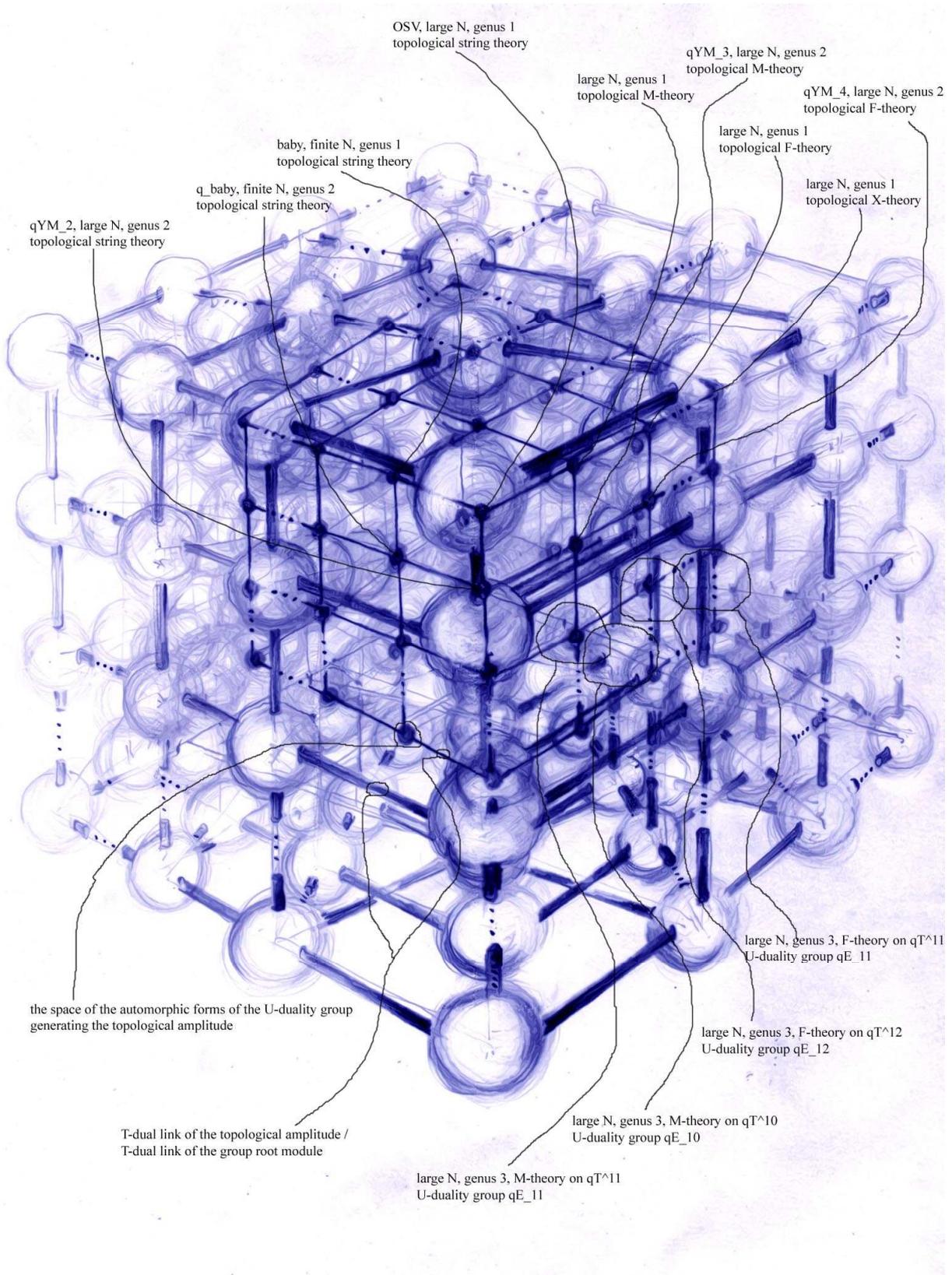

**Definition 2.3** *There is the bijective homomorphism between the inclusion system of the submodules of the graded module of the vacuum torsion and the isomorphism web of the root system of the S-dual modular group*

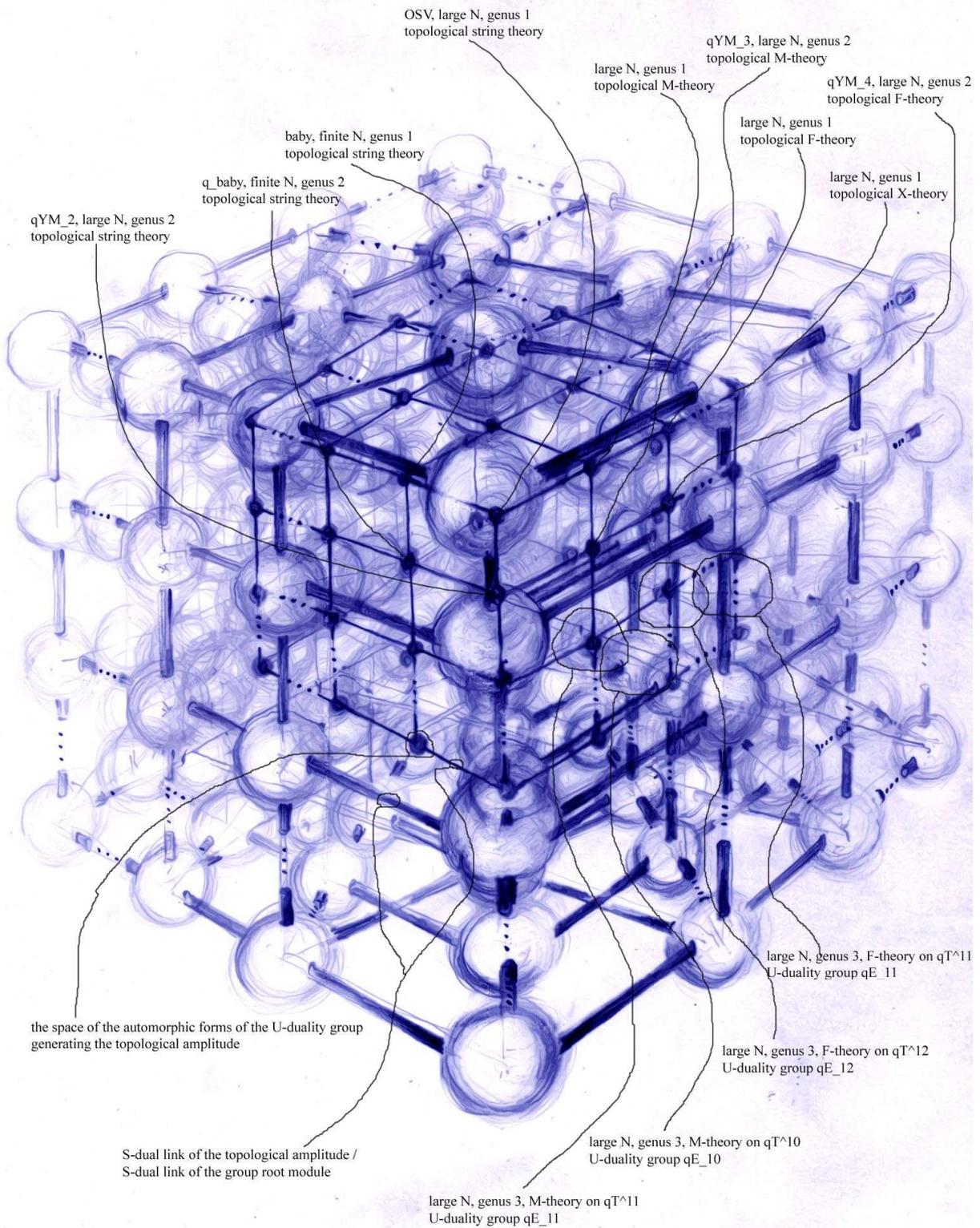

**Definition 2.4** There is the bijective homomorphism between the inclusion system of the submodules of the graded module of the vacuum torsion and the isomorphism web of the root system of the U-dual modular group **M**

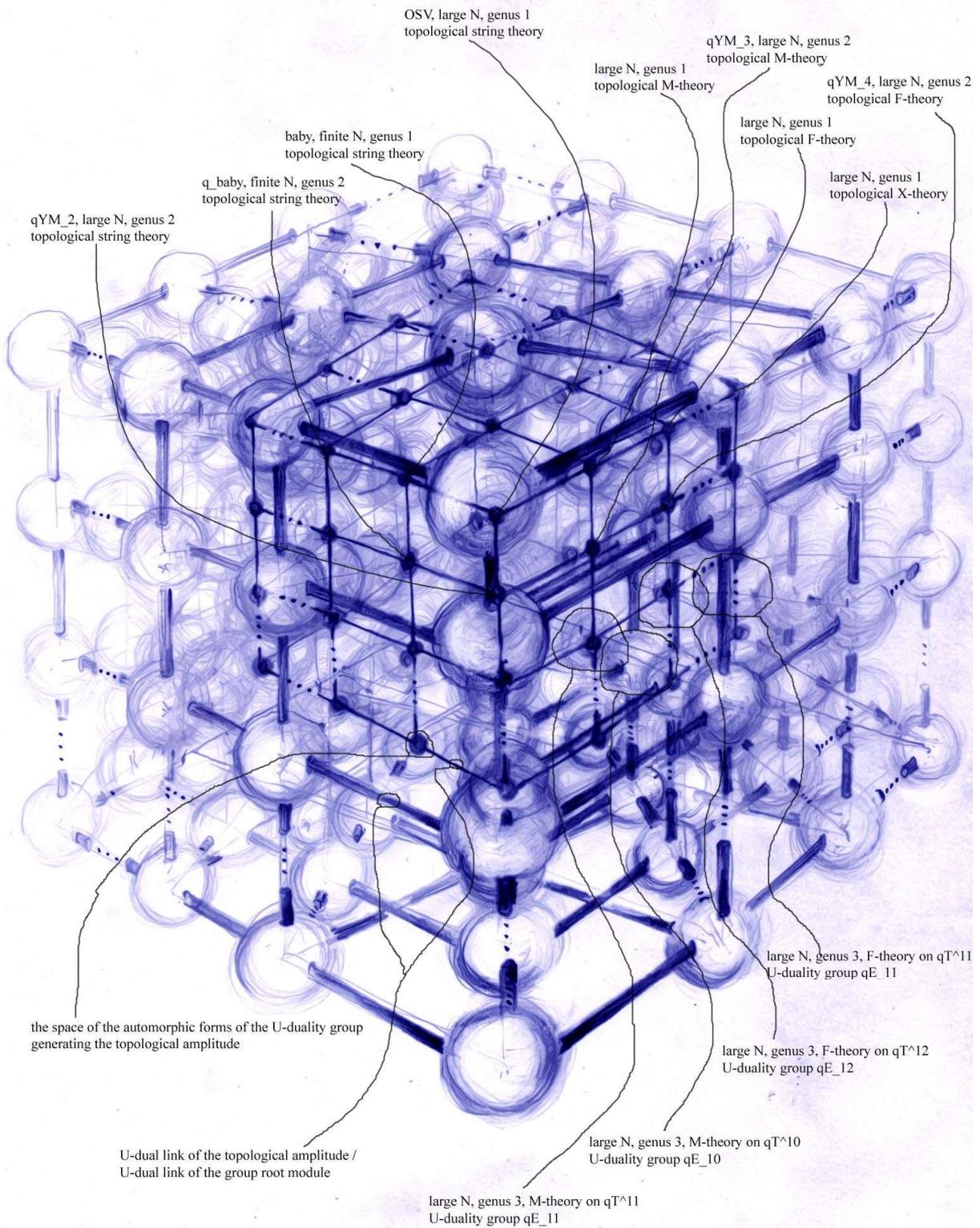

Labels (clockwise from top):
- OSV, large N, genus 1 topological string theory
- qYM_3, large N, genus 2 topological M-theory
- large N, genus 1 topological M-theory
- qYM_4, large N, genus 2 topological F-theory
- large N, genus 1 topological F-theory
- large N, genus 1 topological X-theory
- large N, genus 3, F-theory on qT^11, U-duality group qE_11
- large N, genus 3, F-theory on qT^12, U-duality group qE_12
- large N, genus 3, M-theory on qT^10, U-duality group qE_10
- large N, genus 3, M-theory on qT^11, U-duality group qE_11
- U-dual link of the topological amplitude / U-dual link of the group root module
- the space of the automorphic forms of the U-duality group generating the topological amplitude
- qYM_2, large N, genus 2 topological string theory
- q_baby, finite N, genus 2 topological string theory
- baby, finite N, genus 1 topological string theory

An elegant evidence is Bott periodicity implication for the root lattice of the U-duality group q$E_{27}$ of bosonic M-theory on the quantum deformed 27 dimensional torus q$T^{27}$. The space of the automorphic forms of q$E_{27}$ generate modular form of 19 dimensional bosonic topological M-theory which is just qYM$_{19}$. The root lattice of the U-duality group q$E_{28}$ of bosonic F-theory on q$T^{28}$ is isomorphic to the root lattice of q$E_{27}$

$$\Gamma_{qE_{27}} = \Gamma_{qE_{28}}$$

The space of the automorphic forms of q$E_{28}$ generating modular form of 20 dimensional bosonic topological F-theory, which is qYM$_{20}$, is isomorphic to the modular form of qYM$_4$, which is just topological F-theory.

**Theorem 2.5** *There is a map between the quantum deformed Weyl group of the root space of the intersection cycles and T-dual vacuum torsion.*

*Proof.* There is the isomorphism between the vacuum automorphism symmetries and the T-dual root lattice automorphism symmetries

$$\text{Aut}(\text{vacuum}) = \text{Aut}(\Gamma_{T-dual})$$

Q.E.D.

**Theorem 2.6** *There is a map between the quantum deformed Weyl group of the root space of the intersection cycles and the S-duality cascade of the generalized conifold.*

*Proof.* Coherent superposition of the dilaton condensate, representing the orbifold group G of the generalized conifold, is generated by the S-dual root lattice automorphism group

$$G = \text{Aut}(\Gamma_{S-dual})$$

Q.E.D.

**Theorem 2.7** *There is a map between the quantum deformed Weyl group of the root space of the intersection cycles and the U-dual instanton moduli space permuting the vacuum fluxes.*

*Proof.* U-dual instantons $\zeta_i$ coincide with the simple roots of **M** $\alpha_i$

$$\alpha_i = \zeta_i$$

Q.E.D.

# 3. The vacuum torsion as the automorphism group of the root lattice of M

**Definition 3.1** *Let a lower dimensional brane in an unstable vacuum is a defect of the gauge field of the brane filling the group variety **M**, then a charges of a lower dimensional black brane throat are the excitations of the black brane throat filling **M**.*

N U-dual branes wrapping the U-dual cycle C generate a vacuum flux of the U-dual form Ψ through C

$$N\mathbf{g} = \int_C \Psi$$

We classify the black brane throat fluctuations by the K-theoretical torsion classes of the vacuum variety **M**, so the U-dual fluxes are the quantum deformed Chern characters of the K-theoretical group. Our main observation is that the charges of the U-dual branes wrapping the U-dual cycles are just the roots of **M**.

**Definition 3.2** *The irreducible representation of **M** is coherent superposition of q**O**-valued K-theoretical condensate. The U-dual flux is determined by the U-dual vacuum torsion: the instanton number of the U-dual cycle is defined by the cohomology group of **M**/by the homotopy group of **M**.*

Fact that **M** permutes the instanton numbers of the U-dual cycles is in agreement with our **Theorem 2.7**. Our perspective is twofold: first the basis vectors of the U-dual dilaton φ$_i$ coincide with simple roots of **M**, second Dynkin diagram of **M** represents the intersecting U-dual cycles

```
                o
                |
    o--o--o--o--o--o--o- ... -o

    φ1                    ... φN
```

Our topological result relates the module of the intersecting U-dual cycles determining the vacuum variety potential and the integral of the nonassociative star product of the U-dual forms

$$\#(C_1, C_2, C_3, \ldots, C_N) = \int_M \Psi_1 * \Psi_2 * \Psi_3 * \ldots * \Psi_N$$

**Theorem 3.3** *There is the isomorphism between the orbifold of the generalized conifold throat and the automorphism group of the root lattice of* **M**.

*Proof.* The vacuum multiplets are generated as the orbits of Weyl group **W** of the root space of **M**

$$\mathbf{W} = \text{Aut}(\text{vacuum})$$

Q.E.D.

**Corollary 3.4** *There is the root module homomorphism between the U-dual branes wrapping the U-dual cycles and the orbifold group G of the generalized conifold.*

*Proof.* There is the homomorphism between the S-duality cascade of the generalized conifold and the quantum deformed Weyl reflections of the root space of the vacuum condensate

$$\text{Aut}(\text{vacuum}) = G$$

Q.E.D.

We observe that the root system of **M** generates off-shell arithmetic chaos in brane cosmology. Arithmetic chaos inside Weyl chambers of **M** is just the nilpotent U-fold-valued Wardian flux. We base on the correspondence between the cascading conifold throat and the nilpotent orbits of **M**.

Ward identity of the U-dual dilaton condensate is just the entanglement principle of the vacuum torsion which we observe as the superposition principle of the potential of the moduli space of **M**. The isomorphism between U-dual instanton amplitude and the space of the automorphic forms of **M** gives out mapping of the U-dual instanton-induced vacuum algebraic structure. Thus we base on Feynmanian integral treatment over the U-fold.

**Corollary 3.5** *The generating function of the vacuum torsion $\Xi$ is given by the path integral over* **M**.

*Proof.* The space of the automorphic forms $\Delta$ of **M** generates modular form $\Xi$

$$\Delta = \Xi$$

Q.E.D.

**Theorem 3.6** *There is the quantum deformed Cartan matrix homomorphism between the vacuum torsion and **M**.*

*Proof.* There is the isomorphism between the quantum deformed Cartan matrix $[A_{ij}]_q$ of the vacuum torsion and the S-duality cascade of the generalized conifold which is homomorphic to the monodromy of the group variety **M**

$$[A_{ij}]_q = G$$

Q.E.D.

Next evidence is that

$$[A_{ij}]_q = \int_M \Psi_1 * \Psi_2 * \Psi_3 * \ldots * \Psi_N = \#(C_1, C_2, C_3, \ldots, C_N)$$

**Corollary 3.7** *Ξ is generated by the orbits of the fundamental representation of **M**.*

*Proof.* There is the isomorphism between the module of the intersecting U-dual cycles and Ξ

$$\Xi = [A_{ij}]_q$$

Q.E.D.

We define U-dual modular form Ξ through the one-to-one correspondence between the space of the automorphic forms Δ and Hilbert space H containing the irreducible representations of the vacuum torsion

$$\Delta = H$$

**Definition 3.8** *The microscopic degeneracies of the U-dual black brane throat are the generalized Fourier coefficients of Ξ.*

**Theorem 3.9** *The nilpotent orbits of **M** determine the U-dual black brane amplitude.*

*Proof.* Let β is the weight of **M** and η is dual Coxeter number of **M**, then the coupling constant **g**

$$\mathbf{g} = 2\pi i(\beta\eta)^{-1} = e^{\varphi}$$

and

$$e^{-\mathbf{g}} = q(\mathbf{M})$$

The partition function of the U-dual black brane

$$Z_{U-\text{dual bh}} = \int_M D\mathbf{M}\, |e^{q(\mathbf{M})}|^2$$

Q.E.D.

**Corollary 3.10** *The automorphism group of the root lattice of **M** determines nilpotent U-fold-valued Wardian flux.*

*Proof.* The U-dual black brane amplitude is generated by $\Delta$

$$Z_{U-\text{dual bh}} = \Delta$$

Q.E.D.

Fundamental principle, the superposition principle of the generalized conifold transitions, is just Ward identity of $\Xi$. In the line of **Theorem 3.3** and **Corollary 3.10** we expose $Z_{U\text{-dual bh}}$ as a graceful proof for the quantum deformed Riemann hypothesis. We observe in the fractal attractor of the condensation of the quantum deformed Hagedorn spectrum of the dilaton condensate filling the U-dual black brane throat that a vacua accumulate very near the points of the generalized conifold and formates a fractal flower. Components of the bifurcation diagram for the quantum deformed multi-boundary singularities correspond to the quantum deformed Bernoulli numbers.

As a next evidence for our reasoning we survey the mirror-dual of $\Xi$ as the U-dual instanton amplitude.

**Theorem 3.11** *There is the homomorphism between the root module of the U-dual instanton and the quantum deformed Hagedorn spectrum of the dilaton condensate filling $Z_{U\text{-dual bh}}$.*

*Proof.* There is the isomorphism between the modular form of the U-dual instanton $\zeta$ and the modular form of the U-dual black brane which is determined by the automorphism group of the root lattice of **M**

$$Z_\zeta = Z_{U-\text{dual bh}}$$

Q.E.D.

**Theorem 3.12** *There is the isomorphism between the roots degeneracy of **M** and divergent fluctuations of the quantum deformed Hagedorn spectrum of the dilaton condensate filling $Z_{U\text{-dual bh}}$.*

*Proof.* The basis vectors of the dilaton condensate are $\Gamma_M$-valued vectors

$$\varphi_i = \alpha_i$$

Q.E.D.

**Corollary 3.13** *There is the isomorphism between the vacuum torsion and the uncountable roots degeneracy of M.*

*Proof.* The volume of the U-dual cycle, which is determined by monodromy of the group variety **M**, and the rank of the vacuum torsion, which is determined by **W**, are continuous variables, and

$$\mathbf{W} = [A_{ij}]_q$$

Q.E.D.

It indicates that the conifold throat warping ($=\varphi$) is continuously tuneable parameter and the background flux countinuously varies. Elementary observation is that, thanks to the divergences of the vacuum fluctuations, the propability of the condensation of $q\mathbb{O}$-valued K-theoretic condensate is vanishig.

**Corollary 3.14** *The moduli space potential of **M** and the module of the intersecting cycles living on $\Gamma_M$ are vanishing.*

*Proof.* The uncountable roots degeneracy of **M** coincides with the U-dual black brane degeneracy

$$Z_{U-\text{dual bh}} = \int_\mathbf{M} D\mathbf{M} \, |e^{q(\mathbf{M})}|^2 = 0$$

Q.E.D.

# 4. Conclusion: The symmetry of the vacuum torsion module

Our comprehending must undergo futhermore revolution: Ξ is formulated in Platonian mathematical universe by the uncountable many mathematical consistent formalisms.

**Definition 4.1** *A model of the equipotential surface to the consistency strength of the mathematical formalisms defining the U-dual vacuum torsion which is generated by the nilpotent orbits of* **M**

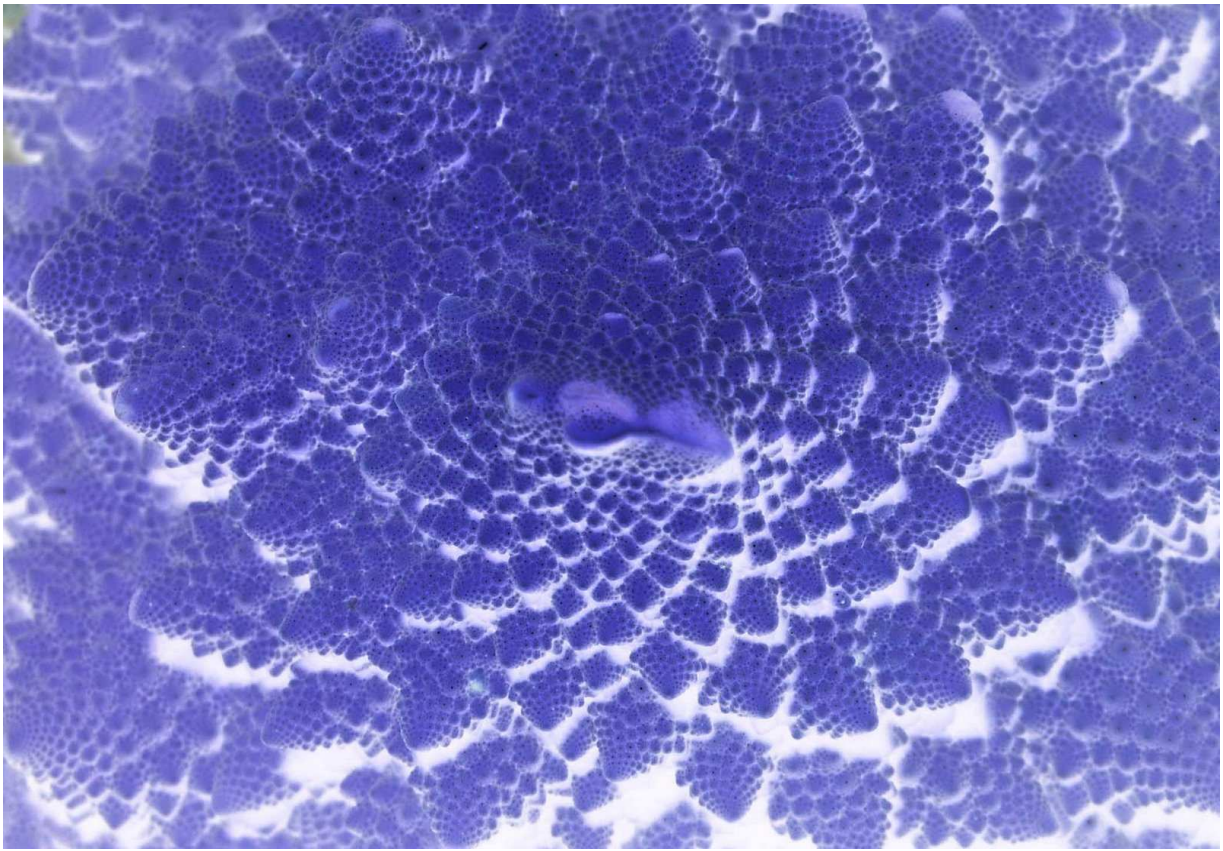

**Definition 4.2** *Generally there is the uncountable many axiomatic systems generating the mathematical consistent formalisms defining the basis module*

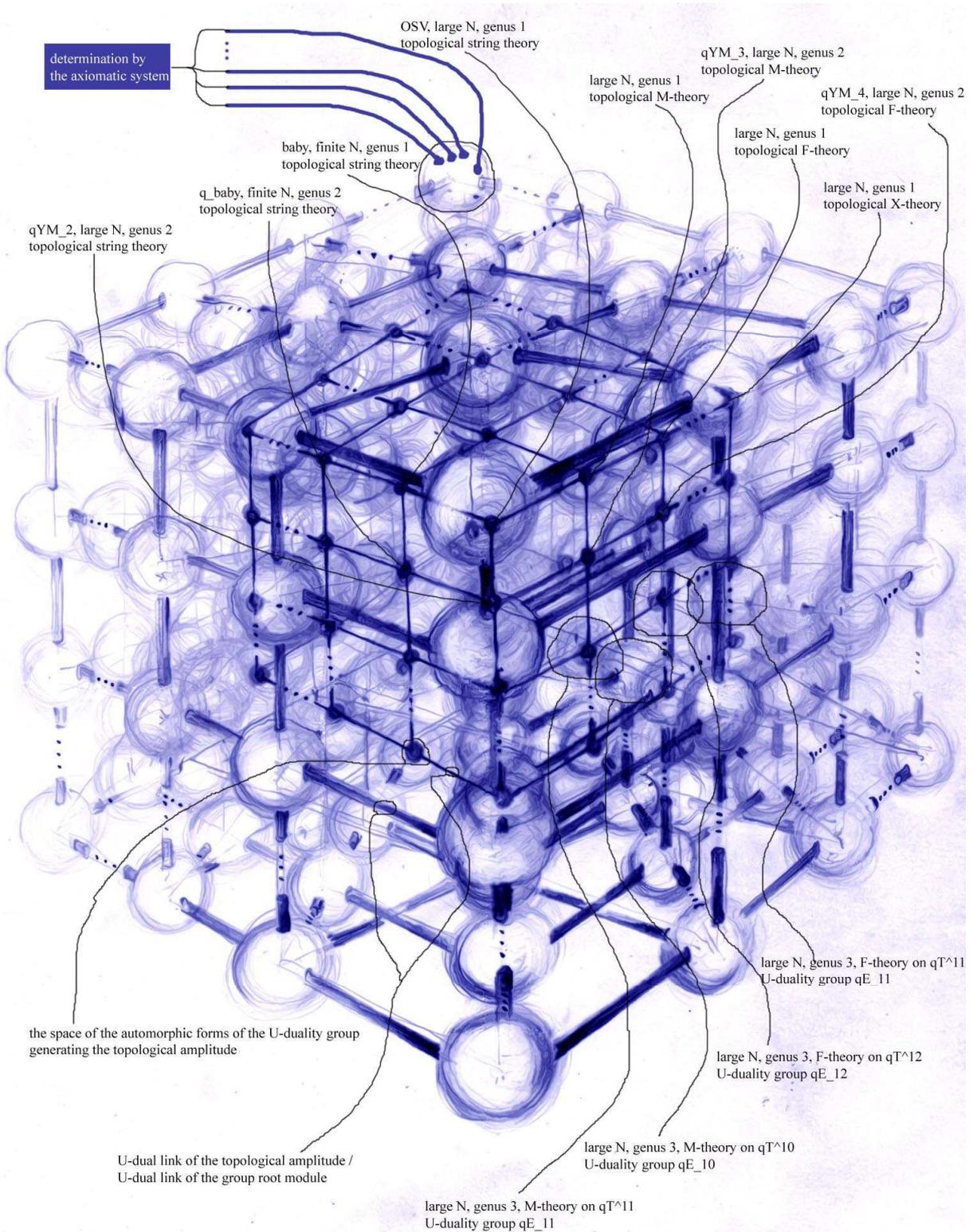

**Notes**

[1] topology_at_atlas.cz

[2] Viz Appendix.

[3] There is $q\mathbf{D_4}$ isomorphism between $q\text{Cliff}_8$ and the endomorphism group $\text{End}(q\mathbb{O} \oplus q\mathbb{O})$.

# 5. Appendix: The quantum deformed group

Let $A_{ij}$ is Cartan matrix of Kac-Moody algebra and q is the quantum parameter. Then the quantum deformed group generate generators $X_i$, $Y_i$ (for the simle roots $\alpha_i$) and $Z_\beta$ (where $\beta$ is an element of the weight lattice)

$$Z_\beta X_i Z_\beta^{-1} = q^{(\beta,\alpha_i)} X_i$$

$$Z_\beta Y_i Z_\beta^{-1} = q^{-(\beta,\alpha_i)} Y_i$$

$$[X_i, Y_i] = \delta_{ij}(Z_i - Z_i^{-1})/(q_i - q_i^{-1})$$

$$\sum_{n=0}^{1-A_{ij}} (-1)^n \, ([1-A_{ij}]_{q_i}!)/([1-A_{ij}-n]_{q_i}! \, [n]_{q_i}!) \, X_i^n X_j X_i^{1-A_{ij}-n} = 0$$

$$\sum_{n=0}^{1-A_{ij}} (-1)^n \, ([1-A_{ij}]_{q_i}!)/([1-A_{ij}-n]_{q_i}! \, [n]_{q_i}!) \, Y_i^n Y_j Y_i^{1-A_{ij}-n} = 0$$

for $i \neq j$, where

$$Z_i = Z_{\alpha_i}$$

and

$$q_i = q^{1/2 \, (\alpha_i, \alpha_i)}$$

The quantum deformed factorial

$$[0]_{q_i}! = 1$$

$$[n]_{q_i}! = \prod_{m=1}^{n} [m]_{q_i}$$

The quantum deformed number

$$[m]_{q_i} = (q_i^m - q_i^{-m})/(q_i - q_i^{-1})$$

The quantum deformed dimension of the group representation

$$\dim_q(R) = \prod_{1 \leq i < j \leq N} [R_i - R_j + j - i]_q / [j - i]_q = \prod_{\alpha^+ > 0} [\alpha^+(\rho + \Lambda_R)]_q / [\alpha^+ \rho]_q$$

(where $\alpha^+$ are the positive roots, $\rho$ is Weyl vector and $\Lambda_R$ is the hightest weight of the representation R.)